\documentclass[twocolumn,10pt]{article}
\usepackage{amsmath}
\usepackage{amsfonts}
\usepackage{stackrel}
\usepackage{graphicx}
\usepackage{caption}
\usepackage{subcaption}
\usepackage{enumerate}
\usepackage{epstopdf}
\newcommand {\R}{\mathbb R}

\newcommand {\Z}{\mathbb Z}
\newcommand {\N}{\mathbb N}
\newcommand {\e}{\mathrm e}
\renewcommand {\i}{\mathrm i}
\newcommand {\Fix}{\mathrm{Fix}}
\renewcommand {\S}{\mathbf{S}}
\renewcommand {\d}{\mathrm{d}}

\newcommand {\kC}{\mathcal{C}}
\newcommand {\kX}{\mathcal{X}}
\newcommand {\kF}{\mathcal{F}}
\newcommand {\Mat}{\text{Mat}}
\newcommand {\w}{\text{w}}
\newcommand {\re}{\text{Re}}

\title{STABILITY OF HOPF BIFURCATIONS IN TIME-DELAYED FULLY-CONNECTED PLL NETWORKS}

\author{Diego P. Ferruzzo Correa
        Escola Polit\'ecnica da\\
        Universidade de S\~ao Paulo\\
	S\~ao Paulo, Brazil\\
    Email: dferruzzo@usp.br\\	
\'Atila M. Bueno
    Universidade Estadual Paulista\\
       J\'ulio de Mesquita Filho\\
        Campus Experimental de Sorocaba\\
        Sorocaba, S\~ao Paulo, Brazil\\
        Email:atila@sorocaba.unesp.br\\%
       Jos\'e R. Castilho Piqueira
    Escola Polit\'ecnica da\\
        Universidade de S\~ao Paulo\\
	S\~ao Paulo, Brazil\\
    Email: piqueira@lac.usp.br
}

\begin{document}

\maketitle    

\begin{abstract}
{\it Dynamics in delayed differential equations (DDEs) is a well studied problem mainly because DDEs arise in models in many areas of science including biology, physiology, population dynamics and engineering. The change of nature in the solutions in the parameter space for a network of Phase-Locked Loop oscillators was studied in \textit{Symmetric bifurcation analysis of synchronous states of time-delayed coupled Phase-Locked Loop oscillators}. Communications in Nonlinear Science and Numerical Simulation, Elsevier BV, Volume 22, Issues 1–3, May 2015, Pages 793-820, where the existence of Hopf bifurcations for both cases, symmetry-preserving and symmetry-breaking synchronization was well stablished. In this work we continue the analysis exploring the stability of small-amplitude periodic solutions emerging near Hopf bifurcations in the Fixed-point subspace, based on the reduction of the infinite-dimensional space onto a two-dimensional center manifold. Numerical simulations are presented  in order to confirm our analitycal results. Although we explore network dynamics of second-order oscillators, results are extendable to higher order nodes.}
\end{abstract}


\section*{INTRODUCTION}
Networks of oscillators have been studied for decades because their models can represent dynamics in a very wide range of fields, as astronomy, biology, neurology, economics, and the stock market,~\cite{GyHori1991, Martins2013, Li2014, Earl2003,Ponzi2000}. Much of the research has focused to understand the influence that changes in the parameter space have over the global dynamics. In that sense, the stability of the synchronization in the network have been explored from many approaches~\cite{Li2014, Yuan2012, Vuellings2014, Giannakopoulos2001}. In this contibution, we study the stability of small-amplitude periodic solutions which emerge near Hopf bifurcations of the symmetry-preserving solutions for a N-node PLL network, using the centre manifold theorem extended to functional differental equations, and the normal form in the center space. This work is mainly based on previous results obtained on the existence of symmetry-preserving and symmetry-breaking Hopf bifurcations in a N-node second-order PLL network, by Ferruzzo et al.~\cite{FerruzzoCorrea2014a}. 

It is important to note that when the lag in the communications between nodes is considered, the ordinary differential equation (ode) system which describe the network dynamics, becomes a delayed differential equation (dde), whose solution lies, in general, in the function space, and its characteristic equation has infinitely many roots. This particular kind of functional differential equations appear in many engineering problems~\cite{Kalmar-Nagy2001, Stone2004, Gilsinn2002}.
 
We consider the Full Phase model introduced in~\cite{FerruzzoCorrea2014a} to analyse stability of small-amplitude periodic orbits near Hopf bifurcations emerging in the parameter space $(\mu,\tau)$  at the Fixed-point subspace, for the non degenerative case ($K>1$). It has been shown that these bifurcations  can occur when a pair of complex conjugate eigenvalues crosses the imaginary axis in iether direction, from the left to the right and from the right to the left. The main approach used for the analysis is the decomposition of the infinite-dimentional space into a 2-dimensional center space spanned by the eigenvectors correponding to the simple imaginary eigenvalues $\lambda=\pm \i\omega$, $\omega>0$, and into an infinite-dimentional space ``orthogonal'' to the first one (the orthogonality condition will be defined below). We will follow closely the theory and procedures presented in~\cite{Kalmar-Nagy2001, Gilsinn2009, Zhao2009, Campbell2009, Guckenheimer1983}.

\section*{The Full-phase model}
The general model for a $N$-node, fully-connected, second-order oscillator network, in terms of the $i$-th node output phase $\phi_i(t)$, is:
\begin{align}
\ddot\phi_i(t)+\mu\dot\phi_i(t)-\mu-\dfrac{K\mu}{N-1}\displaystyle\sum_{\stackrel[j=1]{j\neq i}{}}^Nf(\phi_i,\phi_j)=0,~~~i=1,\ldots,N,
\label{eq:fullphasemodel}
\end{align}
where:
\begin{align}
  f(\phi_i,\phi_j)=\sin(\phi_j(t-\tau)-\phi_i(t))+\sin(\phi_j(t-\tau)+\phi_i(t)).
  \label{eq:f}
\end{align}
$f:\R\times\R\to\R$, $\mu,K,\tau\in\R^+$, and $N\in\N-\{1\}$.

The equilibria $\phi^\pm$, in equation~\eqref{eq:fullphasemodel}, are:
\begin{align}
  \begin{array}{l}
    \phi^+(n)=\dfrac{1}{2}\left(\arcsin\left(-\dfrac{1}{K}\right)+2n\pi\right)\\\\
    \phi^-(n)=\dfrac{1}{2}\left(\pi-\arcsin\left(-\dfrac{1}{K}\right)+2n\pi\right)\\
  \end{array},
  \label{eq:equilibria}
\end{align}
$n\in\mathbb{Z},~K\geq 1$. For our analysis, we consider three main assumptions:
\begin{enumerate}[(a)]
\item\label{aa}The critical eigenvalue $\lambda$ of the linearization of~\eqref{eq:fullphasemodel} at equilibria crosses the imaginary axis with non vanishing velocity, i.e. $\re(\lambda'(\phi^\pm))\neq 0$.
\item\label{ab}The purely imaginary eigenvalue $\lambda=\i\omega$ is simple.
\item\label{ac}The linearization of~\eqref{eq:fullphasemodel} at equilibria, has no eigenvalues of the form $\i k \omega$, $k\in\Z-\{1,-1 \}$.
\end{enumerate}
The Taylor expansion of \eqref{eq:fullphasemodel} at equilibria is:
\begin{align}
  \begin{array}{l}
    \delta\ddot\phi_i+\mu\delta\dot\phi_i-\dfrac{K\mu}{N-1}\displaystyle\sum_{\stackrel[j=1]{j\neq i}{}}^N\sum_{r=1}^\infty\bigg\{\dfrac{1}{r!}\bigg(\delta\phi_i\frac{\partial}{\partial\phi'_i}+\\
    ~~~~~~~~~~~~~\delta\phi_{j\tau}\dfrac{\partial}{\partial\phi'_{j\tau}} \bigg)^rf(\phi_i,\phi_{j\tau}) \bigg\}_{\stackrel[\phi'_{j\tau}=\phi^\pm]{\phi'_i=\phi\pm}{}}=0.
  \end{array}
\end{align}
where $\phi_{j\tau}:=\phi_j(t-\tau)$. Truncate the series up to the third-order term:
\begin{align}
  \begin{array}{l}
    \ddot\phi_i+\mu\dot\phi_i=\dfrac{K\mu}{N-1}\displaystyle\sum_{\stackrel[j=1]{j\neq i}{}}^N\bigg\{\left(\phi_{j\tau}-\phi_i \right)+\left(\phi_{j\tau}+\phi_i \right)\cos 2\phi^\pm\\
    -\dfrac{1}{2}\left(\phi_{j\tau}+\phi_i \right)^2\sin 2\phi^\pm-\dfrac{1}{6}\left[\left(\phi_{j\tau}-\phi_i \right)^3+\left(\phi_{j\tau}+\phi_i \right)^3\cos 2\phi^\pm \right]\bigg\}
  \end{array},
  \label{eq:linear_expansion_form}
\end{align}
$i=1,\ldots,N$, here, for the sake of notation we changed $\delta\phi_i\to\phi_i$. 

The vector field form $\dot x=G(x_t,x;\eta)$, $G:\R^{2N}\times\R^{2N}\times\R^3\to\R^{2N}$, can be obtained by choosing $x^{(i)}_1=\phi_i$ and $x^{(i)}_2=\dot\phi_i$, then, the restriction $G|_i$, or, $\dot x^{(i)}=G^{(i)}(x_t,x^{(i)};\eta)$ gives:
\begin{align}
  \begin{array}{l}
    \dot{x}_1^{(i)}=x_2^{(i)}\\
    \dot{x}_2^{(i)}=-\mu x_2^{(i)}+\dfrac{K\mu}{N-1}\displaystyle\sum_{\stackrel[j=1]{j\neq i}{}}^N\bigg\{(-1+\cos 2\phi^\pm)x^{(i)}_1 \\
    ~~~~~~+(1+\cos 2\phi^\pm)x^{(j)}_{1\tau} -\dfrac{1}{2}\left(x^{(j)}_{1\tau}+x^{(i)}_1 \right)^2\sin 2\phi^\pm\\
    ~~~~~~-\dfrac{1}{6}\left[\left(x^{(j)}_{1\tau}-x^{(i)}_1 \right)^3+\left(x^{(j)}_{1\tau}+x^{(i)}_1 \right)^3\cos 2\phi^\pm \right]\bigg\},\\
i=1,\ldots,N.
  \end{array}
  \label{eq:vector_form}
\end{align}
Following~\cite{Hale1971, Gilsinn2009}, we can represent the dynamics in~\eqref{eq:vector_form} by the abstract differential equation:
\begin{align}
  \frac{d}{dt} x_t(\phi)=A(\eta)x_t(\phi)+\kF(x_t(\phi),\eta).
  \label{eq:abstractODE}
\end{align}
We define $\kX:=\kC([-\tau,0],\R^{2N})$, the Banach space of continuous functions from $[-\tau,0]$ into $\R^{2N}$, equipped with the usual norm $$\|\vartheta\|=\stackrel[-\tau\leq\theta\leq0]{}{\text{sup}}|\vartheta(\theta)|,~~~\vartheta\in\kC([-\tau,0]),$$ $x_t$, in~\eqref{eq:abstractODE}, lies in $\kX$ and satisfies $(T(t)\phi)(\theta)=(x_t(\phi))(\theta)=x(t+\theta)$, where $T(t)$ is a semigroup of family of operators, $\theta\in[-\tau,0]$, and $\eta$ is a vector of parameters. The linear operator $A(\eta)\in\Mat(2N)$ is:
\begin{align}
  (A(\eta)\vartheta)=\left\{\begin{array}{ll}
      \frac{\partial\vartheta}{\partial\theta}(\theta)&,-\tau<\theta\leq 0\\
      A_0(\eta)\vartheta(0)+A_\tau(\eta)\vartheta(-\tau)&,\theta=0
    \end{array}\right.,
  \label{eq:operator_A}
\end{align}
where $A_0(\eta):=\frac{\partial G}{\partial x}\big|_{\phi^\pm}$, $A_\tau(\eta):=\frac{\partial G}{\partial x_\tau}\big|_{\phi^\pm}$
and,
\begin{align}
  \left(\kF(x)\right)(\theta)=\left\{
    \begin{array}{ll}
      \frac{\partial x}{\partial\theta}(\theta)&,-\tau\leq\theta<0\\
      F(x(0),x(-\tau),\eta)&,\theta=0
    \end{array}
  \right.,
  \label{eq:general_F}
\end{align}
$F=(f^{(1)},\ldots,f^{(N)})^T$, $f^{(i)}=(f^{(i)}_1,f^{(i)}_2)$, $f^{(i)}_1=0$, and $f^{(i)}_2=\dfrac{K\mu}{N-1}\displaystyle\sum_{\stackrel[j=1]{j\neq i}{}}^N\bigg\{-\dfrac{1}{2}\left(x^{(j)}_{1\tau}+x^{(i)}_1 \right)^2\sin 2\phi^\pm-\dfrac{1}{6}\left[\left(x^{(j)}_{1\tau}-x^{(i)}_1 \right)^3+\left(x^{(j)}_{1\tau}+x^{(i)}_1 \right)^3\cos 2\phi^\pm \right]\bigg\}.$

In order to build the decomposition of the infinite-dimensional space, we need the adjoint operator associated to the linear part of the linearization and a inner product, via a bilinear form. Associated to the linear part of~\eqref{eq:abstractODE}, the formal adjoint equation is
\begin{align}
  \dfrac{dy}{dt}(t,\eta)=A_0^T(\eta)y(t,\eta)+A_\tau^T(\eta)y(t+\tau,\eta),
  \label{eq:adjoint}
\end{align}
The strongly continuous semigroup $(T^*(t)\psi)(\theta)=(y_t(\psi))(\theta)=y(t+\theta)$, defines the infinitesimal generator:
\begin{align}
  (A^*(\eta)\psi)=\left\{\begin{array}{ll}
      \frac{\partial\psi}{\partial\theta}(\theta)&,0<\theta\leq\tau\\
      A_0(\eta)^T\psi(0)+A_\tau(\eta)^T\psi(\tau)&,\theta=0
    \end{array}\right.,
  \label{eq:adjoint_gen}
\end{align}
such that $\frac{d}{dt}T^*(t)\psi=A^*T^*(t)\psi$, $\psi\in\kX^*:=\kC([0,\tau],\R^{2N})$. The natural inner product has the form~\cite{Hale1993}:
\begin{align*}
  \langle x,y \rangle=\bar{x}^T(0)y(0)+\int_{-\tau}^0\bar{x}^T(s+\tau)A_\tau(\eta)y(s)ds,
\end{align*}
$x\in\kX$ and $y\in\kX^*$; thus, we have~\cite{Gilsinn2009}:
\begin{enumerate}
\item $\lambda$ is an eigenvalue of $A(\eta)$ if and only if $\bar\lambda$ is and eigenvalue of $A^*(\eta)$.
\item If $\varphi_1,\ldots,\varphi_d$ is a basis for the eigenspace of $A(\eta)$ and $\psi_1,\ldots,\psi_d$ is a basis for the eigenspace of $A^*(\eta)$, construct the matrices $\Phi=(\varphi_1,\ldots\varphi_d)$ and  $\Psi=(\psi_1,\ldots,\psi_d)$. Define the bilinear form:
  \begin{align}
    \langle \Psi,\Phi \rangle=I
\label{eq:bilinear_form}
  \end{align}
\end{enumerate}

\section*{The Fixed Point space $\S_N$}
Due to the $\S_N$-symmetry of~\eqref{eq:fullphasemodel} the space where solutions $\phi_i$ lie can be decomposed into the Fixed-point subspace where symmetry-preserving solutions emerge  and a subspace with symmetry-breaking solutions, this was shown in~\cite{FerruzzoCorrea2014a}. We analyze stability of the small-amplitude periodic solutions near Hopf bifurcations in the Fixed point space, these bifurcations satisfy assumptions~\eqref{aa}-\eqref{ac} for $K>1$. In this subspace, equation~\eqref{eq:vector_form} has the form:
\begin{align}
  \begin{array}{l}
    \dot x_1=x_2\\
    \dot x_2=-\mu x_2+K\mu(-1+\cos 2\phi^\pm)x_1\\
            ~~~+K\mu\bigg\{(1+\cos 2\phi^\pm)x_{1\tau}-\dfrac{1}{2}(x_{1\tau}+x_1)^2\sin 2\phi^\pm\\
            ~~~-\dfrac{1}{6}\left[(x_{1\tau}-x_1)^3+(x_{1\tau}+x_1)^3\cos 2\phi^\pm \right]\bigg\},
  \end{array}
\end{align}
then matrices $A_0(\eta)$ and $A_\tau(\eta)$ in~\eqref{eq:adjoint_gen} become:
\begin{align}
A_0(\eta) = \left(
  \begin{array}{cc}
    0&1\\
K\mu(-1+\cos 2\phi^\pm )&-\mu
  \end{array}
\right),
\label{eq:A0}
\end{align}
\begin{align}
A_\tau(\eta)=\left(
  \begin{array}{cc}
    0&0\\
    K\mu(1+\cos 2\phi^\pm)&0
  \end{array}
\right),
  \label{eq:Atau}
\end{align}
and $F$ in~\eqref{eq:general_F} takes the form $F=(f_1~f_2)^T$, with $f_1=0$, and $f_2$:
\begin{align}
\begin{array}{l}
f_2(x_t,\eta)=K\mu\bigg\{-\dfrac{1}{2}(x_{1\tau}+x_1)^2\sin 2\phi^\pm\\
~~~~~~~~~~~~-\dfrac{1}{6}\left[(x_{1\tau}-x_1)^3+(x_{1\tau}+x_1)^3\cos 2\phi^\pm \right]\bigg\}.
\end{array}
  \label{eq:F}
\end{align}
We need the complex eigenfunctions $As(\vartheta)=\i\omega s(\vartheta)$, $A^*n(\theta)=\i\omega n(\theta)$, associated to the critical eigenvalues $\lambda=\i\omega$, and $\bar \lambda=-\i\omega$ with $s(\vartheta) = s_1(\vartheta)+\i s_2(\vartheta)$ and $n(\theta)=n_1(\theta)+\i n_2(\theta)$. These eigenfunctions can be computed solving the boundary value problem $\frac{d}{d\vartheta}s_{1,2}=\mp\omega s_{2,1}(\vartheta)$, and $\frac{d}{d\theta}n_{1,2}=\pm\omega n_{2,1}(\vartheta)$, which, after substituting the operator $A(\eta)$, becomes:
\begin{align}
\begin{array}{rcl}
  A_0(\eta)s_1(0)+A_\tau(\eta)s_1(-\tau)& = &-\omega s_2(0)\\
  A_0(\eta)s_2(0)+A_\tau(\eta)s_2(-\tau)& = &\omega s_1(0)
\end{array}
\label{eq:A0_sol}
\end{align}
and
\begin{align}
\begin{array}{rcl}
  A_0^T(\eta)n_1(0)+A_\tau^T(\eta)n_1(-\tau)& = &\omega n_2(0)\\
  A_0^T(\eta)n_2(0)+A_\tau^T(\eta)n_2(-\tau)& = &-\omega n_1(0)
\end{array},
\label{eq:Atau_sol}
\end{align}
with general solutions:
\begin{align}
\begin{array}{rcl}
s_1(\vartheta) &=& \cos(\omega\vartheta)c_1-\sin(\omega\vartheta)c_2\\
s_2(\vartheta) &=& \sin(\omega\vartheta)c_1+\cos(\omega\vartheta)c_2\\
n_1(\theta)&=& \cos(\omega\theta)d_1-\sin(\omega\theta)d_2\\
n_2(\theta) &=& \sin(\omega\theta)d_1+\cos(\omega\theta)d_2
\end{array}.
\label{eq:gen_sol}
\end{align}
The coefficients $c_1=[c_{11}~c_{12}]^T,~c_2=[c_{21}~c_{22}]^T,~d_1=[d_{11}~d_{12}]^T,~d_2=[d_{21}~d_{22}]^T$ can be obtained by considering the boundary conditions,
\begin{align}
\begin{array}{l}
  \left(\begin{array}{c}
    A_0(\eta)+\cos(\omega\tau)A_\tau(\eta)\\\omega I+\sin(\omega\tau)A_\tau(\eta)
  \end{array}\right)^T\left(
  \begin{array}{c}
    c_1\\c_2
  \end{array}
\right)=0\\\\
 \left( \begin{array}{c}
    A_0^T(\eta)+\cos(\omega\tau)A_\tau^T(\eta)\\
-\omega I-\sin(\omega\tau)A_\tau^T(\eta)
  \end{array}\right)^T\left(
  \begin{array}{c}
    d_1\\d_2
  \end{array}
\right)=0
\end{array},
\end{align}
the ``orthonormality'' condition $\langle s,n\rangle=I$, and setting $c_{11}=1$ and $c_{21}=0$, see ~\cite{Kalmar-Nagy2001, Hale1977} for more details.

It is also possible to decompose the solution $x_t(\vartheta)$ to equation~\eqref{eq:abstractODE} into $x_t(\vartheta)=y_1(t)s_1(\vartheta)+y_2(t)s_2(\vartheta)+\w_t(\vartheta)$, where $y_1$ and $y_2$ lie in the center subspace, such that $y_{1,2}(t)=\langle n_{1,2}(0),x_t(0) \rangle$, and $\w_t$ in the infinite-dimensional component subspace, thus, we have
\begin{align}
  \begin{array}{rcl}
    \dot y_1 &=& \omega y_2 + n_1^T(0)F\\
    \dot y_2 &=& -\omega y_1 + n_2^T(0)F\\
\end{array}
\label{eq:y1y2}
\end{align}
\begin{align}
     \dot \w  = A(\eta)\w_t+\mathcal{F}(x_t,\eta)-n_1^T(0)Fs_1 - n_2^T(0)Fs_2,
\label{eq:trans}
\end{align}
where
\begin{align}
\begin{array}{l}
  \mathcal{F}=\left\{\begin{array}{ll}
      0&,\vartheta\in[-\tau,0)\\
      F(y_1(t)s_1(0)+y_2(t)s_2(0)+\w(t)(0))&,\vartheta=0.
    \end{array}\right.
\end{array}
\label{eq:FF-F}
\end{align}

\subsection*{The center manifold}
Following~\cite{Hassard1981, Kalmar-Nagy2001, Zhao2009}, we know that $\w$ can be approximated by the second-order expansion:
\begin{align}
 \w(y_1,y_2)(\vartheta)=\dfrac{1}{2}(h_1(\vartheta)y_1^2+2h_2(\vartheta)y_1y_2+h_3(\vartheta)y_2^2),
\label{eq:center_manifold}
\end{align}
thus, by differentiating and substituting equation~\eqref{eq:trans} keeping up to second order terms, we obtain:
\begin{align}
  \dot\w = -\omega h_2y_1^2+\omega(h_1-h_3)y_1y_2 + \omega h_2 y_2^2  + O(y^3),
\label{eq:dot_w}
\end{align}
and from equation~\eqref{eq:trans}, 
\begin{align}
  \dfrac{d\w}{dt}=A(\eta)\w_t+\mathcal{F}(\w+y_1s_1+y_2s_2)-(d_{12}s_1+d_{22}s_2)f_2.
\label{eq:dwdt}
\end{align}
From the definition of $A(\eta)$, equivalent to~\eqref{eq:adjoint_gen}, we see that
\begin{align}
  A(\eta)\w=\left\{\begin{array}{ll}
      \frac{1}{2}(\dot h_1 y_1^2+2\dot h_2 y_1y_2 + \dot h_3y_2^2)&,\vartheta\in[-\tau,0)\\
      A_0(\eta)\w(0)+A_\tau(\eta)\w(-\tau)&,\vartheta=0,
    \end{array}\right .
\label{eq:Aw}
\end{align}
then, from equation~\eqref{eq:center_manifold},~\eqref{eq:dot_w},~\eqref{eq:dwdt}, and~\eqref{eq:Aw}, we can obtain the unknown coefficients $h_1,~h_2$, and $h_3$, solving:
\begin{align}
  \begin{array}{rcl}
    \dot h_1 &=& 2( -\omega h_2 +f_2^{20}(d_{12}s_1(\vartheta)+d_{22}s_2(\vartheta))),\\
    \dot h_2 &=& \omega(h_1-h_3)+ f_2^{11}(d_{12}s_1(\vartheta)+d_{22}s_2(\vartheta)),\\
    \dot h_3 &=& 2(\omega h_2 + f_2^{02}(d_{12}s_1(\vartheta)+d_{22}s_2(\vartheta))),\\
  \end{array}
\label{eq:h1h2h3a}
\end{align}
and,
\begin{align}
  \begin{array}{l}
A_0(\eta)h_1(0)+ A_\tau(\eta) h_1(-\tau)=\\
~~~~~~~~~2(-\omega h_2(0)+ f_2^{20}(d_{12}s_1(0)+d_{22}s_2(0))),\\\\
A_0(\eta)h_2(0)+A_\tau h_2(-\tau)=\\
~~~\omega (h_1(0)-h_3(0)) + f_2^{11}(d_{12}s_1(0)+d_{22}s_2(0))),\\\\
A_0(\eta)h_3(0)+ A_\tau(\eta) h_3(-\tau)=\\
~~~~~~~~~2(\omega h_2(0)+ f_2^{02}(d_{12}s_1(0)+d_{22}s_2(0))),
  \end{array}
\label{eq:h1h2h3}
\end{align}
where $f^{20} = \dfrac{1}{2}\dfrac{\partial^2 f}{\partial y_1^2}\bigg|_0$, $f^{11} = \dfrac{\partial^2 f}{\partial y_1\partial y_2}\bigg|_0$, and $f^{02} = \dfrac{1}{2}\dfrac{\partial^2 f}{\partial y_2^2}\bigg|_0$.

Equation~\eqref{eq:h1h2h3a} is written as the inhomogeneus differential equation:
\begin{align}
  \frac{dh}{d\vartheta}= Ch+p\cos(\omega\vartheta)+q\sin(\omega\vartheta)
\label{eq:inhomogeneus_h}
\end{align}
where
\begin{align*}
\begin{array}{l}
  h:= \left(
    \begin{array}{c}
      h_1\\h_2\\h_3
    \end{array}
\right),~~~C:= \omega\left(
  \begin{array}{ccc}
    0&-2 I &0\\
    I&0&-I\\
    0&2 I&0
  \end{array}
\right)_{6\times 6}\\\\
p:=\left(
  \begin{array}{c}
    f_2^{20}p_0\\f_2^{11}p_0\\f_2^{02}p_0
  \end{array}
\right),~~~q:=\left(
  \begin{array}{c}
    f_2^{20}q_0\\f_2^{11}q_0\\f_2^{02}q_0
  \end{array}
\right),\\\\
p_0:=\left(
  \begin{array}{c}
    d_{12}\\c_{22}d_{22}
  \end{array}
\right),~~~q_0:=\left(
  \begin{array}{c}
    d_{22}\\-c_{22}d_{12}
  \end{array}
\right),
\end{array}
\end{align*}
with general solution:
\begin{align}
h(\vartheta)=\e^{C\vartheta}K + M\cos(\omega\vartheta) + N\sin(\omega\vartheta).
  \label{eq:gen_sol}
\end{align}
After substituting the general solution  into~\eqref{eq:inhomogeneus_h} we solve for $M$ and $N$, and then from the boundary value problem we solving for $K$,
 \begin{align}
\left(
  \begin{array}{cc}
    C&-\omega I\\
   \omega I& C
  \end{array}
\right)\left(
  \begin{array}{c}
    M\\N
  \end{array}
\right)=-\left(
  \begin{array}{c}
    p\\q
  \end{array}
\right)
   \label{eq:MN}
 \end{align}
 \begin{align}
   Ph(0)+Qh(-\tau)=p-r,
\label{eq:K_sol}
 \end{align}
where
\begin{align}
\begin{array}{l}
  P:=\left(
    \begin{array}{ccc}
      A_0&0&0\\
      0&A_0&0\\
      0&0&A_0
    \end{array}
\right)-C,\\\\
Q:=\left(
    \begin{array}{ccc}
      A_\tau&0&0\\
      0&A_\tau&0\\
      0&0&A_\tau
    \end{array}
\right),
\end{array}
\end{align}
and $r:=\left(
  \begin{array}{cccccc}
    0&f_2^{20}&0&f_2^{11}&0&f_2^{02}
  \end{array}
\right)^T$.

The expressions for $\w_1(0)$ and $\w_1(-\tau)$, necessary in~\eqref{eq:FF-F}, are:
\begin{align}
\begin{array}{ll}
  \w_1(0)&=\dfrac{1}{2}\bigg( (M_1+K_1)y_1^2 + 2(M_3+K_3)y_1y_2+(M_5+K_5)y_2^2\bigg),\\\\
 \w_1(-\tau)&=\dfrac{1}{2}\bigg((\e^{-C\tau}K|_1+M_1\cos(\omega\tau)-N_1\sin(\omega\tau))y_1^2\\
&~~~+2(\e^{-C\tau}K|_3+M_3\cos(\omega\tau)-N_3\sin(\omega\tau))y_1y_2\\
&~~~+(\e^{-C\tau}K|_5+M_5\cos(\omega\tau)-N_5\sin(\omega\tau))y_2^2\bigg),
\end{array}
\label{eq:W0Wtau}
\end{align}
note that we only need $\w_1(\vartheta)$ since the nonlinear function in~\eqref{eq:F} only depends on $x_1$; then by substituting~\eqref{eq:W0Wtau} into~\eqref{eq:y1y2}, we obtain: 
\begin{align}
\begin{array}{rr}
\dot y_1 =& \omega y_2 + g_1(y_1,y_2;\eta)\\
\dot y_2 =& -\omega y_1 + g_2(y_1,y_2;\eta)
\end{array},
\end{align}
or
\begin{align}
  \begin{array}{l}
    \dot y_1=\omega y_2 + a_{20}y_1^2 +a_{11}y_1y_2 + a_{02}y_2^2 + a_{30}y_1^3\\
~~~~~~+a_{21}y_1²y_2+a_{12}y_1y_2^2+ a_{03}y_2^3,\\
    \dot y_2=-\omega y_1 + b_{20}y_1^2 +b_{11}y_1y_2 + b_{02}y_2^2 + b_{30}y_1^3\\
~~~~~~+b_{21}y_1²y_2+b_{12}y_1y_2^2+ b_{03}y_2^3.
  \end{array}
\label{eq:doty1y2}
\end{align}
In \cite{Guckenheimer1983} is computed the coefficient $a$, which determines stability of the normal form~\eqref{eq:doty1y2},  
\begin{align}
\begin{array}{ll}
  a&=\dfrac{1}{16}\left[g^{03}_2+g^{21}_2+g^{12}_1+g^{30}_1 \right]+\dfrac{1}{16\omega}\bigg[g^{11}_2\left(g^{02}_2+g^{20}_2 \right)\\
   &~~~-g^{11}_1\left(g^{02}_1+g^{20}_1 \right)- g^{02}_2g^{02}_1+g^{20}_2g^{20}_1\bigg],
\end{array}
\label{eq:a}
\end{align}
where $g^{ij}_r=\dfrac{\partial^{i+j}}{\partial^iy_1\partial^jy_2}g_r(0,0)$. Periodic orbits near Hopf bifurcation at the critical eigenvalue $\lambda=\i\omega$, will be stable if $a<0$ and unstable if $a>0$.

\section*{Numerical Results}
\begin{figure}[!htb]
\centering
  \includegraphics[width=0.5\textwidth]{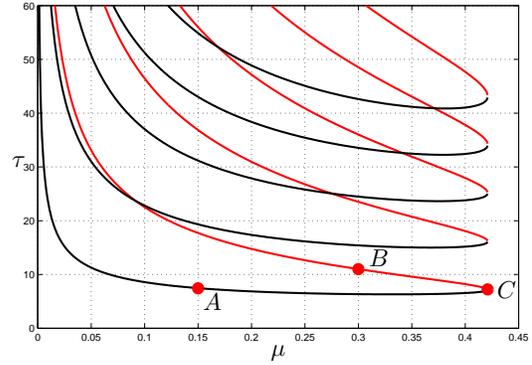}
  \caption{\MakeUppercase{Symmetry-preserving bifurcation curves in $\Fix(\S_N)$ for $K=1.05$. In black, Hopf bifurcations with $\re(\lambda')>0$ and, in red, Hopf bifurcations with $\re(\lambda')<0$.}}
  \label{fig:symm_bif_Fix}
\end{figure}

We reproduced some of the computations for the Hopf bifurcations in the Fixed point space for the case $K>1$ presented in~\cite{FerruzzoCorrea2014a}, because we will compute stability for these bifurcation curves using results obtained in the previous section. In figure~\ref{fig:symm_bif_Fix} (part of figure 10, in~\cite{FerruzzoCorrea2014a}) are shown the symmetry-preserving bifurcations curves in the parameter space $(\mu,\tau)$ for $K=1.05$, for both cases: bifurcations with $\re(\lambda')>0$ in black color, and with $\re(\lambda')<0$ in red color; we  also choose three testing point for numerical simulation $A=(\mu,\tau)=(0.15,7.46)$, $B=(0.3,11)$, and $C=(0.421,7.10)$.

\begin{figure}[!htb]
\centering
  \includegraphics[width=0.5\textwidth]{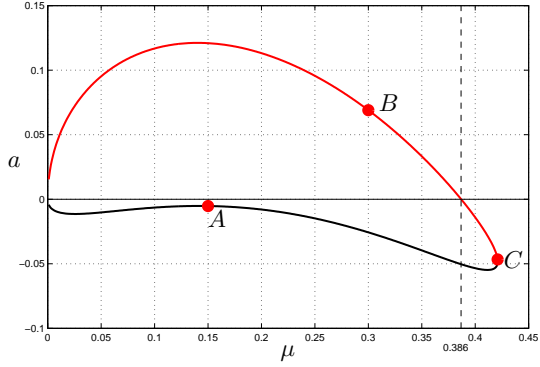}
\caption{COEFFICIENT $a$ COMPUTED USING EQ.\eqref{eq:a}, DETERMINING STABILITY OF HOPF BIFURCATIONS IN $\Fix(\S_N)$ FOR $K=1.05$, SEE FIGURE~\ref{fig:symm_bif_Fix}.}
\label{fig:a}
\end{figure}

In figure~\ref{fig:a} is shown the coefficient $a$ computed using equation~\eqref{eq:a}, in the parameter space $(\mu,\tau)$, for $K=1.05$, related to the Hopf bifurcations curves shown in figure~\ref{fig:symm_bif_Fix}. The  black curve corresponds to stability of periodic orbits near Hopf bifurcations with $\re(\lambda')>0$ (black curves in figure~\ref{fig:symm_bif_Fix}), as we can see, these periodic solutions are all stable ($a<0$). The red curve corresponds to stability of periodic orbits near Hopf bifurcations with $\re(\lambda')<0$ (red curves in figure~\ref{fig:symm_bif_Fix}), these periodic orbits are unstable for $\mu<\mu^*(K)\approx 0.386$, and stable for $\mu>\mu^*$. Are also shown point $A$, $B$ anc $C$. At points $A$ and $C$, small amplitude periodic orbits are stable, whilst at point $B$, they are unstable.  

In order to confirm our results, we computed branches of periodic solutions near the Hopf bifurcations points $A$, $B$, and $C$, using DDE-BIFTOOL~\cite{Engelborghs2001, Engelborghs2002}, along with the Floquet multipliers for a specific periodic solution chosen in the branch.

\begin{figure}
\centering
 \begin{minipage}[b]{.5\textwidth}
  \centering
  \includegraphics[width=0.8\textwidth]{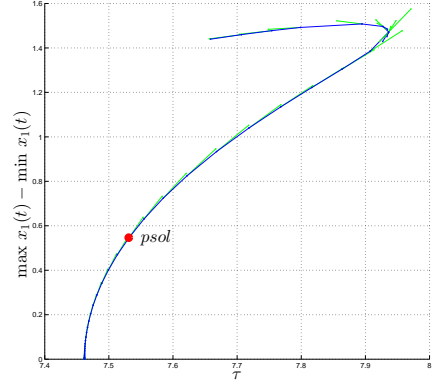}
  \subcaption{}
 \end{minipage}
 \begin{minipage}[!h]{.23\textwidth}
  \includegraphics[width=\textwidth]{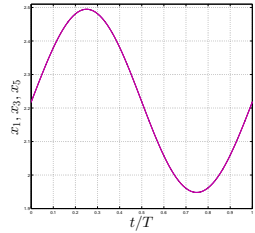}
  \subcaption{}%
 \end{minipage}
 \begin{minipage}[!h]{.23\textwidth}
  \includegraphics[width=\textwidth]{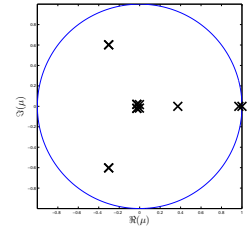}
  \subcaption{}
 \end{minipage}
\caption{(a) BRANCH OF PERIODIC SOLUTIONS EMERGING FROM POINT $A=(\mu,\tau)=(0.15,7.46197)$. (b) PERIODIC SOLUTION PROFILE AT $\mu=0.15$, $\tau=7.5315$, $T=12.0364~seg$, (POINT $psol$). (c) FLOQUET MULTIPLIERS FOR THE PERIODIC SOLUTION $psol$.}
\label{fig:pointA}
\end{figure}

A branch of periodic solutions with small amplitude, emerging from the Hopf bifurcation point $A=(\mu,\tau)=(0.15,7.46)$ is shown in figure~\ref{fig:pointA}-(a). In~\ref{fig:pointA}-(b) it is shown the periodic solution profile $psol$ at $\tau=7.5315$. The Floquet multipliers related to $psol$ are shown in figure~\ref{fig:pointA}-(c). It is clear that this periodic solution is stable, since there is no Floquet multiplier outside the unity circle.

\begin{figure}
\centering
 \begin{minipage}[b]{.5\textwidth}
  \centering
  \includegraphics[width=0.8\textwidth]{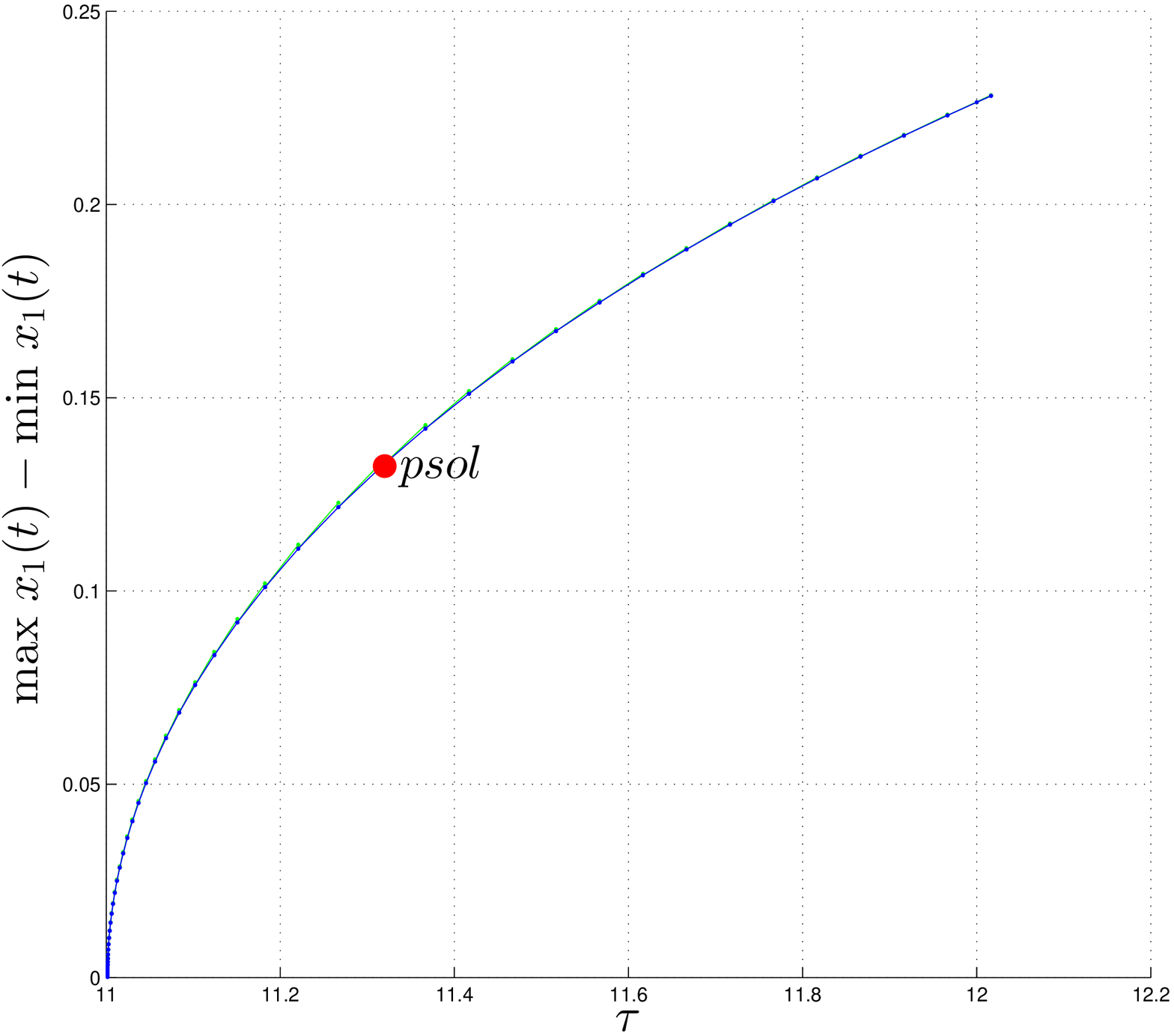}
  \subcaption{}
 \end{minipage}
 \begin{minipage}[!h]{.23\textwidth}
  \includegraphics[width=\textwidth]{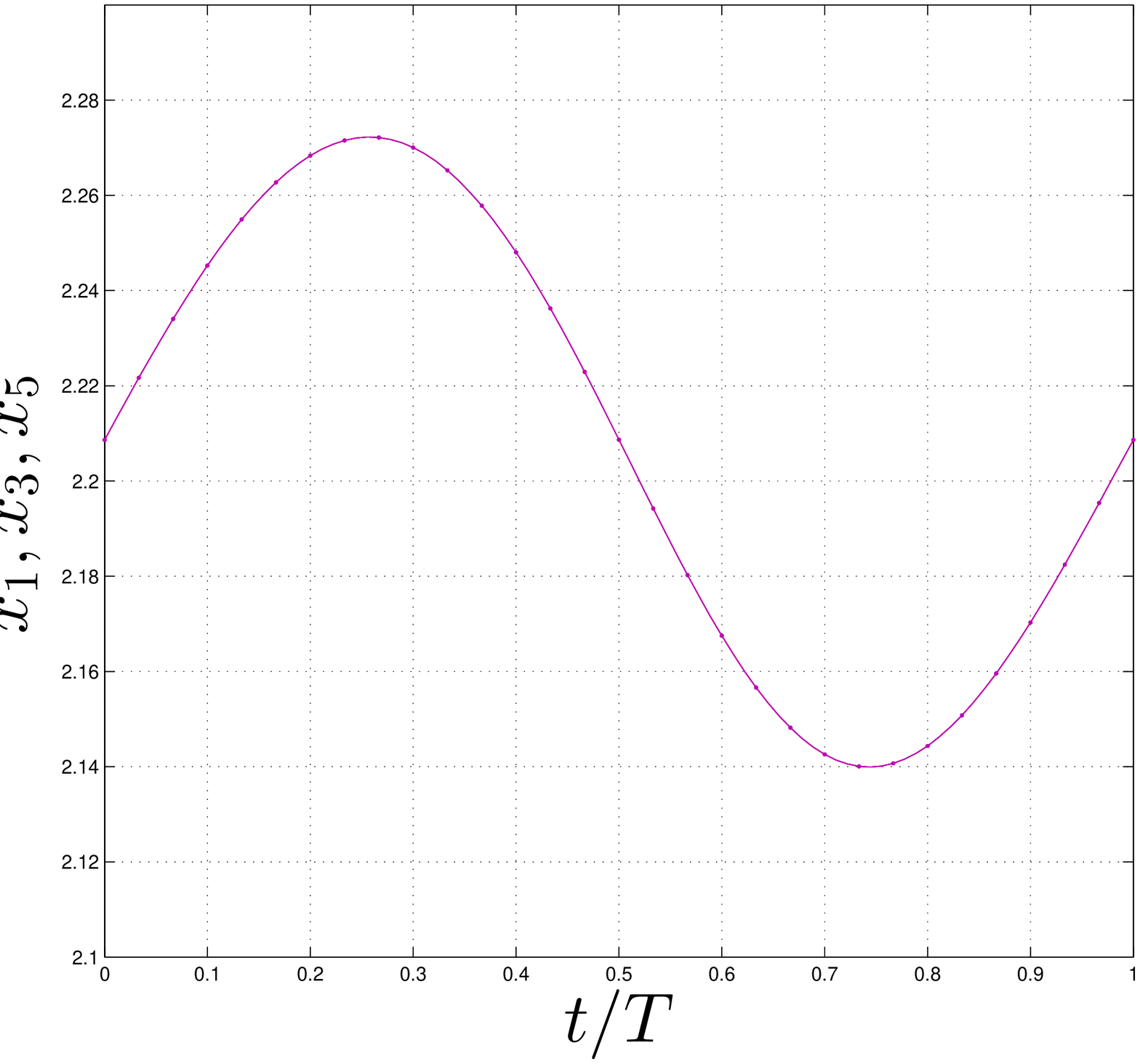}
  \subcaption{}%
 \end{minipage}
 \begin{minipage}[!h]{.23\textwidth}
  \includegraphics[width=\textwidth]{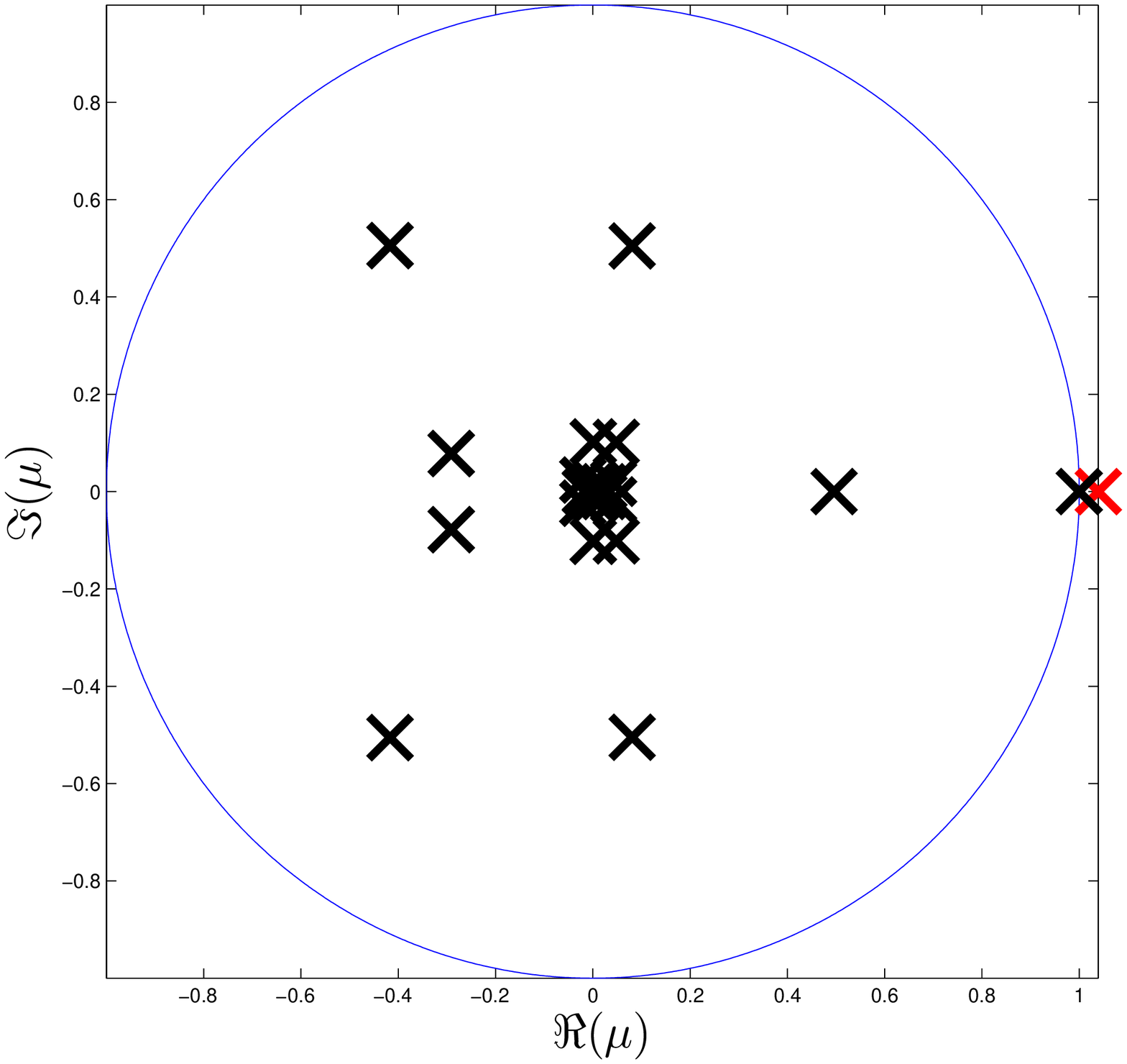}
  \subcaption{}
 \end{minipage}
\caption{(a) BRANCH OF PERIODIC SOLUTIONS EMERGING FROM POINT $B=(\mu,\tau)=(0.3,11.001518)$. (b) PERIODIC SOLUTION PROFILE AT $\mu=0.3$, $\tau=11.3744$, $T=12.8506~seg$, (POINT $psol$). (c) FLOQUET MULTIPLIERS FOR THE PERIODIC SOLUTION $psol$.}
\label{fig:pointB}
\end{figure}

For the point $B=(\mu,\tau)=(0.3,11)$, the branch of periodic solutions is shown in figure~\ref{fig:pointB}-(a), in figure~\ref{fig:pointB}-(b), it is shown the profile for the periodic solution $psol$ chosen at $\tau=11.3744$, these solution is unstable, because there is a Floquet multiplier outside the unity circle, see figure~\ref{fig:pointB}-(c).

\begin{figure}
\centering
 \begin{minipage}[b]{.5\textwidth}
  \centering
  \includegraphics[width=0.8\textwidth]{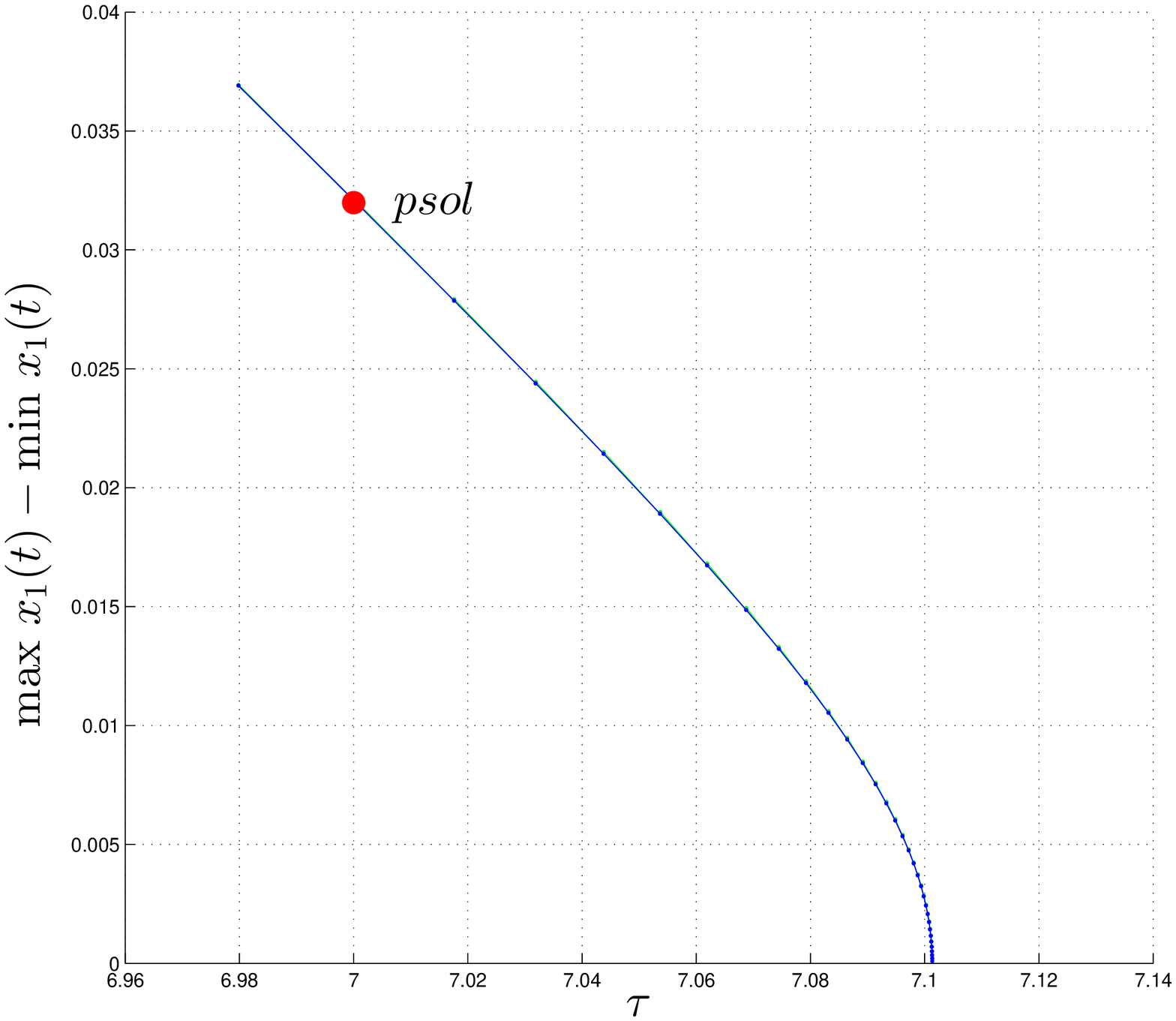}
  \subcaption{}
 \end{minipage}
 \begin{minipage}[!h]{.23\textwidth}
  \includegraphics[width=\textwidth]{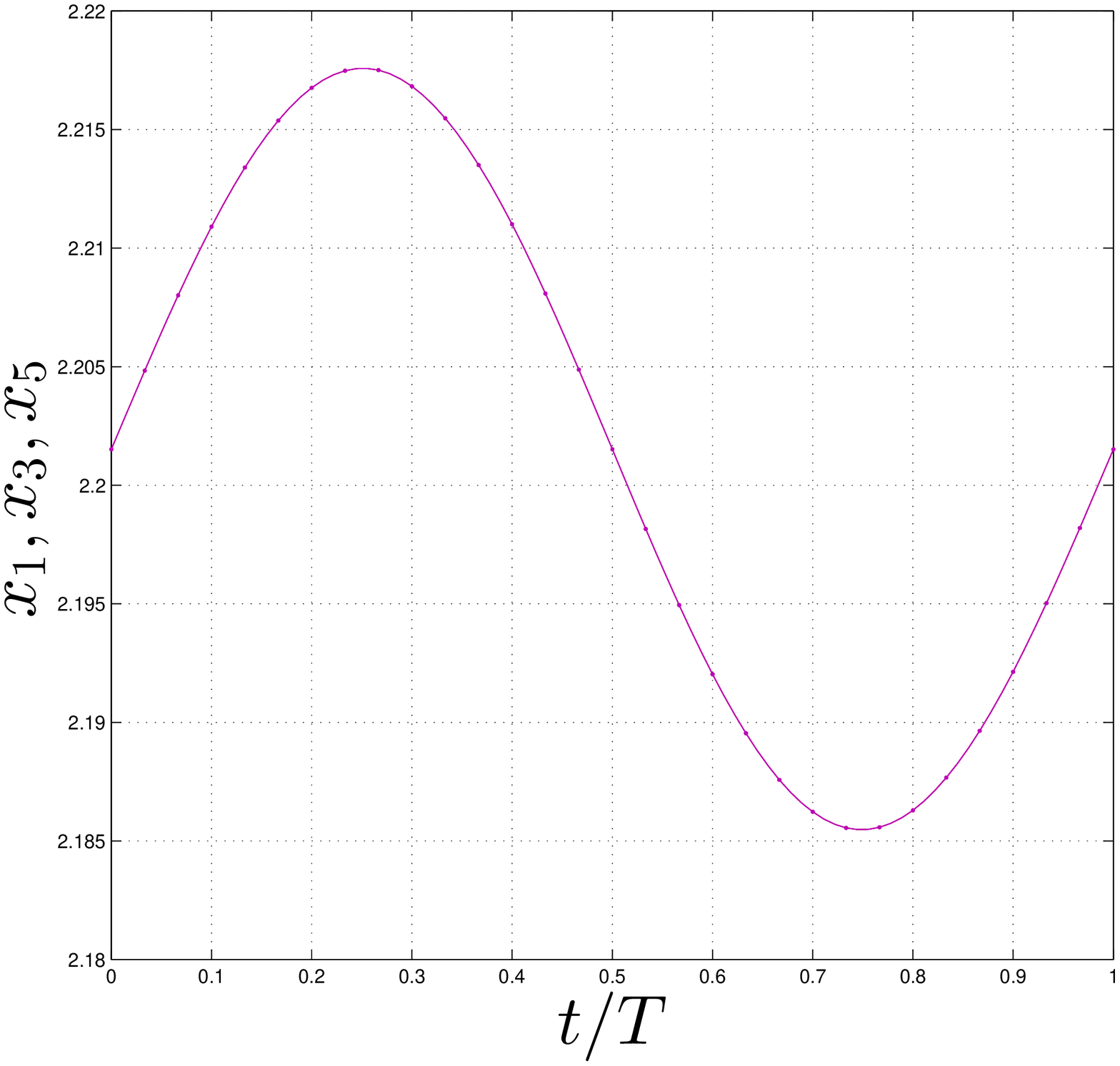}
  \subcaption{}%
 \end{minipage}
 \begin{minipage}[!h]{.23\textwidth}
  \includegraphics[width=\textwidth]{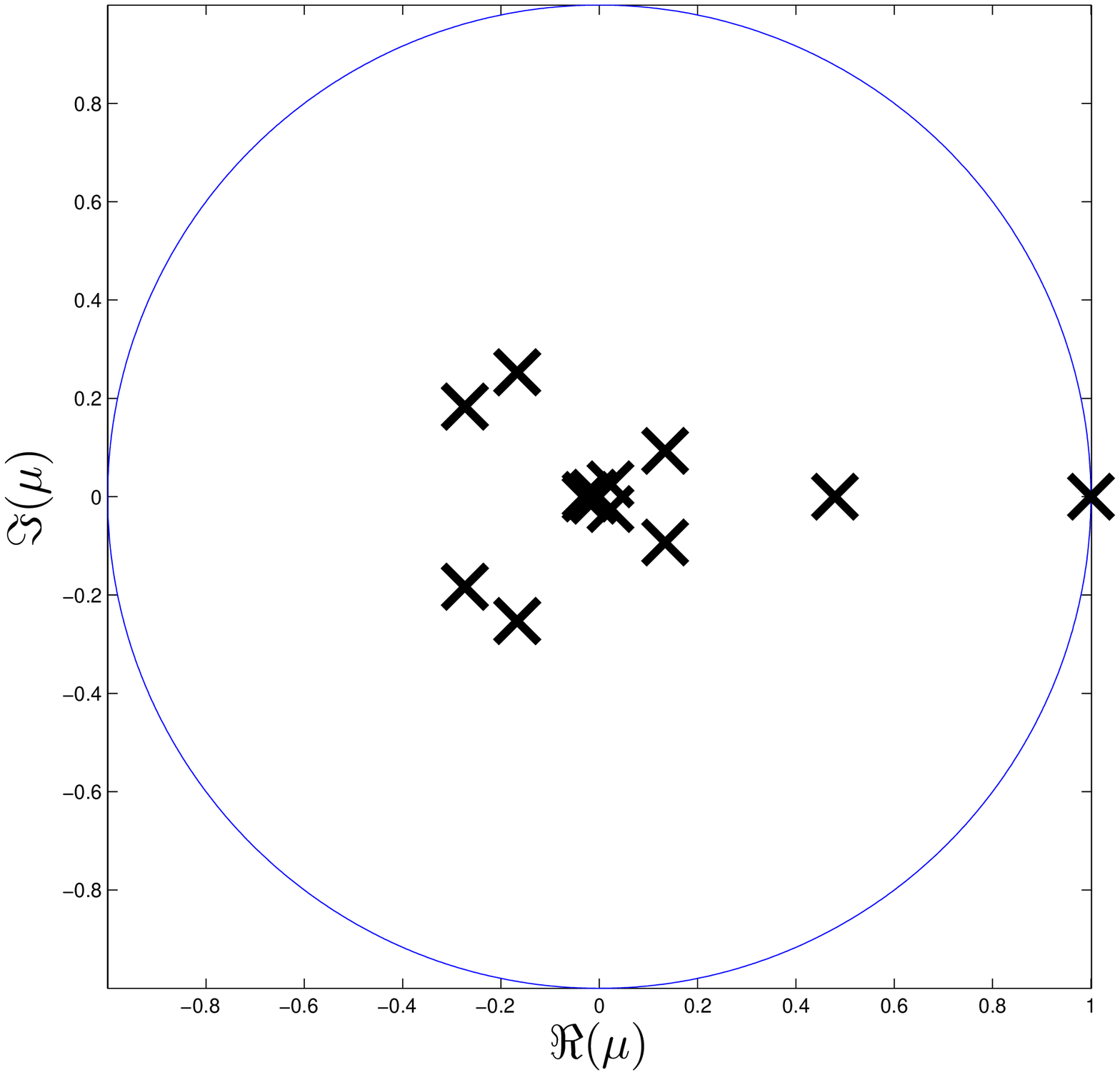}
  \subcaption{}
 \end{minipage}
\caption{(a) BRANCH OF PERIODIC SOLUTIONS EMERGING FROM POINT $C=(\mu,\tau)=(0.421,7.101329)$. (b) PERIODIC SOLUTION PROFILE AT $\mu=0.421$, $\tau=7.00$, $T=8.8704~seg$, (POINT $psol$). (c) FLOQUET MULTIPLIERS FOR THE PERIODIC SOLUTION $psol$.}
\label{fig:pointC}
\end{figure}

Finally, the branch of periodic solutions near the Hopf bifurcation point $C=(\mu,\tau)=(0.421,7.10)$ is shown in figure~\ref{fig:pointC}-(a). The periodic solution choosen in the branch is at $\tau=7.00$, its profile is shown infigure~\ref{fig:pointC}-(b). All the Floquet multipliers shown in figure~\ref{fig:pointC}-(c) are within the unity circle, therefore the solution is stable.

%
\section*{Conclusions}
The reduction of the inifinite-dimensional space onto the center manifold in normal form, was applied to the Fixed point space for the Full-phase model in order to analyse the stability of small-amplitude periodic orbits near simple Hopf bifurcations, in both cases, for $\re(\lambda')>0$ and $\re(\lambda')<0$, we found that in the first case  periodic orbits which are stable ($a<0$) can emerge, and, in the other case, unstable ($a>0$) periodic orbits can emerge for $\mu<\mu^*(K)$, and stable periodic orbits for $\mu>\mu^*(K)$. The numerics show that the analytical results are correct.

Although, we computed the coefficient $a$ for a specific value of $K$, the procedure shown is valid for all the parameter space where simple Hopf bifurcations appear. 

Finally, it is important to spotlight some points for further research: First, what is the nature of the solutions at the special point $\mu=\mu^*(K)$, at which the coefficient $a$ changes sign. Second, analyze stability of the degenerate Hopf bifurcations at the  Fixed point space for $K=1$, which are codimension 2, pure imaginary eigenvalue and zero eigenvalue; and third, the stability of the symmetry-breaking  degenerate Hopf  bifurcations which have multiplicity $N-1$.
\section*{acknowledgment}
We would like to thank the Escola Polit\'ecnica da Universidade de S\~ao Paulo and FAPESP for their support.
%
% ...

\bibliographystyle{plain}
\bibliography{/home/diego/Dropbox/bibliography/phD_bibliography}

\begin{thebibliography}{10}

\bibitem{Campbell2009}
SueAnn Campbell.
\newblock Calculating centre manifolds for delay differential equations using
  maple.
\newblock In {\em Delay Differential Equations}, pages 1--24. Springer US,
  2009.

\bibitem{Earl2003}
Matthew Earl and Steven Strogatz.
\newblock Synchronization in oscillator networks with delayed coupling: A
  stability criterion.
\newblock {\em Phys. Rev. E}, 67(3), Mar 2003.

\bibitem{Engelborghs2002}
K.~Engelborghs, T.~Luzyanina, and D.~Roose.
\newblock Numerical bifurcation analysis of delay differential equations using
  {DDE-BIFTOOL}.
\newblock {\em ACM Trans. Math. Softw.}, 28(1):1--21, March 2002.

\bibitem{Engelborghs2001}
Koen Engelborghs, Tatyana Luzyanina, and Giovanni Samaey.
\newblock {DDE-BIFTOOL} v. 2.00: a {Matlab} package for bifurcation analysis of
  delay differential equations.
\newblock In {\em Numerical Analysis and Applied Mathematics Section}.
  Department of Computer Science, K.U.Leuven, Leuven, Belgium, October 2001.

\bibitem{FerruzzoCorrea2014a}
Diego~Paolo Ferruzzo~Correa, Claudia Wulff, and Jos\'e~Roberto
  Castilho~Piqueira.
\newblock Symmetric bifurcation analysis of synchronous states of time-delayed
  coupled phase-locked loop oscillators.
\newblock {\em Communications in Nonlinear Science and Numerical Simulation},
  22(1–3):793 -- 820, Aug 2015.

\bibitem{Giannakopoulos2001}
Fotios Giannakopoulos and Andreas Zapp.
\newblock Bifurcations in a planar system of differential delay equations
  modeling neural activity.
\newblock {\em Physica D: Nonlinear Phenomena}, 159(3-4):215 -- 232, 2001.

\bibitem{Gilsinn2009}
David~E. Gilsinn.
\newblock Bifurcations, center manifolds, and periodic solutions.
\newblock In {\em Delay differential equations}, pages 155--202. Springer, New
  York, 2009.

\bibitem{Gilsinn2002}
DE~Gilsinn.
\newblock {Estimating critical Hopf bifurcation parameters for a second-order
  delay differential equation with application to machine tool chatter}.
\newblock {\em {NONLINEAR DYNAMICS}}, {30}({2}):{103--154}, {OCT} {2002}.

\bibitem{Guckenheimer1983}
John Guckenheimer and Philip Holmes.
\newblock {\em Nonlinear oscillations, dynamical systems, and bifurcations of
  vector fields}, volume~42.
\newblock New York Springer Verlag, 1983.

\bibitem{GyHori1991}
I.~Gy{\H{o}}ri and G.~Ladas.
\newblock {\em Oscillation theory of delay differential equations}.
\newblock Oxford Mathematical Monographs. The Clarendon Press Oxford University
  Press, New York, 1991.
\newblock With applications, Oxford Science Publications.

\bibitem{Hale1971}
J.~K. Hale.
\newblock {\em Functional differential equations}.
\newblock Springer-Verlag, New York, 1971.

\bibitem{Hale1993}
J.~K. Hale and S.~M. Verduyn~Lunel.
\newblock {\em Introduction to functional differential equations}.
\newblock Springer-Verlag, London, 1993.

\bibitem{Hale1977}
Jack~K. Hale.
\newblock {\em Theory of Functional Differential Equations (Applied
  Mathematical Sciences)}.
\newblock Springer, 1977.

\bibitem{Hassard1981}
Brian~D. Hassard, Nicholas~D. Kazarinoff, and Yieh~Hei Wan.
\newblock {\em Theory and applications of {H}opf bifurcation}, volume~41 of
  {\em London Mathematical Society Lecture Note Series}.
\newblock Cambridge University Press, Cambridge, 1981.

\bibitem{Kalmar-Nagy2001}
Tam\'as Kalm\'ar-Nagy, G\'abor St\'ep\'an, and Francis~C. Moon.
\newblock Subcritical {H}opf {B}ifurcation in the {D}elay {E}quation {M}odel
  for {M}achine {T}ool {V}ibrations.
\newblock {\em Nonlinear Dynamics}, 26(2):121--142, 2001.

\bibitem{Li2014}
Wenxue Li, Hongwei Yang, Liang Wen, and Ke~Wang.
\newblock Global exponential stability for coupled retarded systems on
  networks: A graph-theoretic approach.
\newblock {\em Communications in Nonlinear Science and Numerical Simulation},
  19(6):1651--1660, Jun 2014.

\bibitem{Martins2013}
A.~Martins and L.H.A. Monteiro.
\newblock Frequency transitions in synchronized neural networks.
\newblock {\em Communications in Nonlinear Science and Numerical Simulation},
  18(7):1786--1791, Jul 2013.

\bibitem{Ponzi2000}
A.~Ponzi and Y.~Aizawa.
\newblock Self-organized criticality and partial synchronization in an evolving
  network.
\newblock {\em Chaos, Solitons \& Fractals}, 11(7):1077 -- 1086, 2000.

\bibitem{Stone2004}
E~Stone and SA~Campbell.
\newblock {Stability and bifurcation analysis of a nonlinear DDE model for
  drilling}.
\newblock {\em {JOURNAL OF NONLINEAR SCIENCE}}, {14}({1}):{27--57}, {JAN-FEB}
  {2004}.

\bibitem{Vuellings2014}
Andrea Vuellings, Eckehard Schoell, and Benjamin Lindner.
\newblock {Spectra of delay-coupled heterogeneous noisy nonlinear oscillators}.
\newblock {\em {EUROPEAN PHYSICAL JOURNAL B}}, {87}({2}), {FEB 3} {2014}.

\bibitem{Yuan2012}
Yuan Yuan and Xiao-Qiang Zhao.
\newblock Global stability for non-monotone delay equations (with application
  to a model of blood cell production).
\newblock {\em Journal of Differential Equations}, 252(3):2189--2209, Feb 2012.

\bibitem{Zhao2009}
Siming Zhao and Tam\'ass Kalm\'ar-Nagy.
\newblock {C}enter {M}anifold {A}nalysis of the {D}elayed {L}ienard {E}quation.
\newblock In {\em Delay Differential Equations}, pages 1--17. Springer US,
  2009.

\end{thebibliography}
%
\end{document}